\documentstyle[12pt]{article}
\catcode`\@=11
\@addtoreset{equation}{section}

\catcode`\@=12
\newtheorem{Theorem}{Theorem}[section]
\newtheorem{Definition}{Definition}[section]
\newtheorem{Proposition}{Proposition}[section]
\newtheorem{Lemma}{Lemma}[section]

\title{A Proof On 
Arnold's Chord Conjecture On Cotangent Bundles
\thanks{Project 19871044 Supported by NSF}}

\author{Renyi Ma \\
Department of Mathematics \\
Tsinghua University \\
Beijing, 100084\\
People's Republic of China\\
rma@math.tsinghua.edu.cn}

\date { }

\begin{document}
\textwidth=165mm
\textheight=210mm
\parindent=8mm
\frenchspacing
\maketitle

\begin{abstract}
In this article, we prove that there 
exists at least one chord which is characteristic 
of Reeb vector field connecting a given Legendre submanifold 
in a contact manifolds of induced type in 
the cotangent bundles of any smooth open manifolds which 
confirms the Arnold conjecture in cotangent 
bundles. 
\end{abstract}


\noindent{\bf 2000 MR Subject Classification} 32Q65, 53D35, 53D12

\section{Introduction and results}

Let $\Sigma$ be a smooth closed oriented manifold of dimension
$2n-1$. A contact form on $\Sigma$ is a $1-$form such that
$\lambda \wedge (d\lambda )^{n-1}$ is a volume form on $\Sigma$.
Associated to $\lambda$ there is 
the so-called Reeb vectorfield $X_\lambda $ defined
by
$$i_X\lambda   \equiv 1, \ \ i_Xd\lambda  \equiv 0.$$

\vskip 3pt 

Concerning the dynamics of Reeb flow, 
there is a well-known conjecture raised by Arnold in \cite{ar} 
which concerned the Reeb orbit and Legendre submanifold in 
a contact manifold. If $(\Sigma ,\lambda )$ is a contact manifold 
with contact form $\lambda $ of dimension $2n-1$, then a Legendre 
submanifold is a submanifold ${\cal L}$ of $\Sigma $, which is 
$(n-1)$dimensional and everywhere tangent to the contact 
structure $\ker \lambda $. Then a characteristic 
chord for $(\lambda ,{{\cal {L}}})$ is a smooth path
$x:[0,T]\to \Sigma ,T>0$
with
$\dot x(t)=X_{\lambda }(x(t)) \ for \ t\in(0,T)$,$x(0),x(T)\in {\cal {L}}$ 
Arnold raised the following conjecture: 

\vskip 3pt 

{\bf Conjecture1}(see\cite{ar}). Let $\lambda _0$ be the standard tight 
contact form $\lambda _0={{1}\over {2}}(x_1dy_1-y_1dx_1+x_2dy_2-y_2dx_2)$
on the three sphere 
$S^3=\{ (x_1,y_1,x_2,y_2)\in R^4|x_1^2+y_1^2+x_2^2+y_2^2=1\}.$ 
If $f:S^3\to (0,\infty )$ is a smooth function 
and ${\cal {L}}$ is a Legendre 
knot in $S^3$, then 
there is a characteristic chord for $(f\lambda _0,{\cal {L}})$. 

\vskip 3pt 

This conjecture 
was completely solved  
in \cite{ma,mo}. 

\vskip 3pt 

In this paper we improve the Gromov's 
proof on that there exists at least one 
intersection point for the weakly 
exact Lagrangian submanifold under the weakly 
Lagrangian isotopy \cite[$2.3.B_{3-4}$]{gro} to prove:

\begin{Theorem}
Let $(\Sigma ,\lambda )$ be a contact manifold with 
contact form $\lambda $ of induced type or Weinstein type in the 
cotangent bundles of any open smooth manifold with symplectic 
form $\sum _{i=1}^ndp_i\wedge dq_i$ induced by 
Liouville form $\alpha =\sum _{i=1}^np_idq_i$, i.e., 
there exists a transversal vector field $Z$ to $\Sigma $ 
such that $L_Z\omega =\omega $, $\lambda =i_Z\omega $.   
Let $X_{\lambda } $ its Reeb vector 
field and ${\cal {L}}$ a closed Legendre submanifold. Then 
$p_idq_i-\lambda  $ defines an element $[p_idq_i-\lambda ]\in H^1(\Sigma )$. 
If $[p_idq_i-\lambda ]=0$, then 
there exists at least one characteristic chord for 
$(X_\lambda ,{\cal {L}})$. 
\end{Theorem}

\vskip 3pt

{\bf Sketch of proofs}: We work in the framework 
as in \cite{gro,ma}. In Section 2, we study the 
linear Cauchy-Riemann operator and sketch some basic 
properties. In section 3, first we construct a 
Lagrangian submanifold $W$ under the assumption that 
there does not exist Reeb chord connecting the Legendre 
submanifold $\cal {L}$; second, we study 
the space ${\cal D}(V,W)$
consisting of contractible disks in manifold $V$ with boundary 
in Lagrangian submanifold $W$ and construct a Fredholm
section of tangent bundle of  ${\cal  D }(V,W)$.  
In section 4, following \cite{gro}, we construct a non-proper 
Fredholm 
section by using a 
special anti-holomorphic section as in \cite{gro,ma}.  
In section 5, we transform the non-homogenious 
Cauchy-Riemann equations to 
$J-$holomorphic curves. 
In section 6, we finish the proof of Theorem 1.1. as in \cite{gro}.

\section{Linear Fredholm Theory}

For $100<k<\infty $ consider the Hilbert space
$V_k$ consisting of all maps $u\in H^{k,2}(D, C^n)$,
such that $u(z)\in R^n\subset C^n$ for almost all $z\in
\partial D$. $L_{k-1}$ denotes the usual Hilbert $L_{k-1}-$space
$H_{k-1}(D, C^n)$. We define an operator
$\bar \partial :V_p\mapsto L_p$ by
\begin{equation}
\bar \partial u=u_s+iu_t
\end{equation}
where the coordinates on $D$ are
$(s,t)=s+it$, $D=\{ z||z|\leq 1\}  $.
The following result is well known(see\cite{al,wen}).
\begin{Proposition}
$\bar \partial :V_p\mapsto L_p $ is a surjective real
linear Fredholm operator of index $n$. The kernel
consists of the constant real valued maps.
\end{Proposition}
Let $(C^n, \sigma =-Im(\cdot ,\cdot ))$ be the standard
symplectic space. We consider a real $n-$dimensional plane
$R^n\subset C^n$. It is called Lagrangian if the
skew-scalar product of any two vectors of $R^n$ equals zero.
For example, the plane $\{(p,q)|p=0\}$ and $\{(p,q)|q=0\}$ are Lagrangian
subspaces. The manifold of all (nonoriented) Lagrangian subspaces of
$R^{2n}$ is called the Lagrangian-Grassmanian $\Lambda (n)$.
One can prove that the fundamental group of
$\Lambda (n)$ is free cyclic, i.e.
$\pi _1(\Lambda (n))=Z$. Next assume
$(\Gamma (z))_{z\in \partial D}$ is a smooth map
associating to a point $z\in \partial D$ a Lagrangian
subspace $\Gamma (z)$ of $C^n$, i.e.
$(\Gamma (z))_{z\in \partial D}$ defines a smooth curve
$\alpha $ in the Lagrangian-Grassmanian manifold $\Lambda (n)$.
Since $\pi _1(\Lambda (n))=Z$, one have
$[\alpha ]=ke$, we call integer $k$ the Maslov index
of curve $\alpha $ and denote it by $m(\Gamma )$, see(\cite{ar}).

Now let $z:S^1\mapsto R^n\subset C^n$ be a smooth curve. 
Then
it defines a constant loop $\alpha $ in Lagrangian-Grassmanian
manifold $\Lambda (n)$. This loop defines
the Maslov index $m(\alpha )$ of the map
$z$ which is easily seen to be zero.

  Now Let $(V,\omega )$ be a symplectic manifold and 
$W\subset V$ a Lagrangian submanifold. Let 
$u:(D, \partial D)\to (V,W)$ be a smooth map homotopic to constant map 
$u_0:(D,\partial D)\to p\in W$. 
Then $u^*TV$ is a symplectic vector bundle and 
$(u|_{\partial D})^*TW$ be a Lagrangian subbundle in 
$u^*TV$. Since $u$ is homotopic to 
$u_0$ by $h(t,z)$ with $h(0,\cdot )=u_0$ and $h(1,\cdot )=u$, 
we can take 
a trivialization of $h^*TV$ as 
$$\Phi (h^*TV)=[0,1]\times D\times C^n$$ 
and 
$$\Phi (h|_{\partial D})^*TW\subset [0,1]\times S^1\times C^n.$$
Let 
$$\pi _2: [0,1]\times D\times C^n\to C^n$$
then 
$$\bar h: (t,z)\in [0,1]\times S^1\to \pi _2\Phi (h|_{\partial D})^*TW|(t,z)\in \Lambda (n).$$
\begin{Lemma}
Let $u: (D^2,\partial D^2) \rightarrow (V,W)$ be a $C^k-$map $(k\geq 1)$ 
as above. Then,
$$m(\bar u)=0$$
\end{Lemma}
Proof.  Since the homotopy $h(t,z)$ induces 
a homotopy $\bar h$ in Lagrangian-Grassmanian
manifold. Note that $m(\bar h(0, \cdot ))=0$.
By the homotopy invariance of Maslov index,
we know that $m(\bar u)=0$.

\vskip 5pt 

   Consider the partial differential equation
\begin{eqnarray}
\bar \partial u+A(z)u=0  \ on \ D  \\
u(z)\in \Gamma (z) R^n\ for \ z\in \partial D \\
\Gamma (z)\in GL(2n,R)\cap Sp(2n)\\
m(\Gamma )=0 \ \ \ \ \ \ \ \ 
\end{eqnarray}

For $100<k<\infty $ consider the Banach space $\bar V_k $
consisting of all maps $u\in H^{k,2}(D, C^n)$ such 
that  $u(z)\in \Gamma (z)$ for almost all $z\in
\partial D$. Let $L_{k-1}$ the usual Hilbert 
space $H_{k-1}(D,C^n)$. 

We define an operator $P$:
$\bar V_{k}\rightarrow L_{k-1}$ by
\begin{equation}
P(u)=\bar \partial u+Au
\end{equation}
where $D$ as in (2.1).
\begin{Proposition}
$\bar \partial : \bar V_k \rightarrow L_{k-1}$
is a real linear Fredholm operator of index n.
\end{Proposition}
Proof: see \cite{al,gro,wen}.

\section{Nonlinear Fredholm Theory}

\subsection{Construction of Lagrangian Submanifold}

Let $M$ be an open manifold and $(T^*M,p_idq_i)$ be 
the cotangent bundle of open manifold with the 
Liouville form $p_idq_i$. Since 
$M$ is open, there exists a function $g:M\to R$ without 
critical point. The translation by 
$tTdg$ along the fibre gives a hamiltonnian isotopy 
of $T^*M$:

\begin{equation}
h^T_t(q,p)=(q,p+tTdg(q))
\end{equation}

\begin{equation}
h^{T*}_t(p_idq_i)=p_idq_i+tTdg.
\end{equation}

\begin{Lemma}
For any given compact set $K\subset T^*M$, there exists 
$T=T_K$ such that $h^T_1(K)\cap K=\emptyset $.
\end{Lemma}
Proof. Similar to \cite{gro,ls}

\vskip 3pt 

Let $\Sigma \subset T^*M$ be a closed hypersurface, if there 
exists a vector field $V$ defined in the neighbourhood $U$ 
of $\Sigma $ transversal to $\Sigma $ 
such that $L_V\omega =\omega $, here 
$\omega =dp_i\wedge dq_i$ is a standard symplectic form  
on $T^*M$ induced by the Liouville form $p_idq_i$, we call 
$\Sigma $ the contact manifold of induced type in $T^*M$ 
with the induced contact form $\lambda =i_V\omega $. 

\vskip 3pt 

Let $(\Sigma ,\lambda )$ be a contact manifold of induced   
type or Weinstein's type in $T^*M$ with contact form 
$\lambda $ and $X$ its Reeb vector field, then 
$X$ integrates to a Reeb flow $\eta _s$ for $s\in R^1$.

By using the transversal vector field $V$, one 
can identify the neighbourhood $U$ of $\Sigma $ 
foliated by flow $f_t$ of $V$ and $\Sigma $, 
i.e., $U=\cup _tf_t(\Sigma )$ with the 
neighbourhood of $\{ 0\}\times \Sigma $ in the symplectization
$R\times \Sigma $ by the exact symplectic transformation(see\cite{ma}).

Consider the form $d(e^a\lambda )$ 
at the point $(a,x)$ on the manifold 
$(R\times \Sigma )$, then one can check 
that $d(e^a\lambda )$ is a symplectic 
form on $R\times \Sigma $. Moreover 
One can check that 
\begin{eqnarray}
&&i_X(e^a\lambda )=e^a \\ 
&&i_X(d(e^a\lambda ))=-de^a 
\end{eqnarray}
So, the symplectization of Reeb vector field $X$ is the 
Hamilton vector field of $e^a$ with 
respect to the symplectic form $d(e^a\lambda )$. 
Therefore the Reeb flow lifts to the Hamilton flow 
$h_s$ on $R\times \Sigma $(see\cite{ag,eg}).

\vskip 3pt 

Let ${\cal L}$ be a closed Legendre submanifold 
in $(\Sigma ,\lambda )$, i.e.,  
there exists a smooth embedding $Q:{\cal L}\to \Sigma $ such that 
$Q^*\lambda |_{\cal L}=0$. Let
$$(V',\omega ')=(T^*M,dp_i\wedge dq_i)$$
and 
\begin{equation}
W'={\cal L}\times R, \ \ W'_s={\cal L}\times \{ s\} \label{eq:3.w}
\end{equation}
define 
\begin{eqnarray}
&&G':W'\to V'  \cr 
&&G'(w')=G'(l,s)=(0,\eta _s(Q(l))) \label{eq:3.ww}
\end{eqnarray}
\begin{Lemma}
There does not exist any Reeb chord connecting Legendre  
submanifold $\cal {L}$
in $(\Sigma ,\lambda )$ if and only if  
$G'(W'(s))\cap G'(W'(s'))$ is empty for $s\ne s'$.
\end{Lemma}
Proof. Obvious. 
\begin{Lemma}
If there does not exist any Reeb chord for $(X_\lambda ,{\cal {L}})$
in $(\Sigma ,\lambda )$ then 
there exists a smooth embedding 
$G':W'\to V'$ with $G'(l,s)=(0,\eta _s(Q(l)))$
such that
\begin{equation}
G'_K:{\cal L}\times (-K, K)\to V'  \label{eq:3.www}
\end{equation}
is a regular open Lagrangian embedding for any finite positive $K$.
\end{Lemma}
Proof. One check 
\begin{equation}
{G'}^*(d(e^a\lambda ))=\eta (\cdot ,\cdot )^*d\lambda 
=(\eta _s^*d\lambda +i_Xd\lambda \wedge ds)=0
\end{equation}
This implies that ${G}'$ is a Lagrangian embedding, this proves 
Lemma3.3. 

\vskip 3pt 

Note that $d\lambda =dp_i\wedge dq_i $ on 
$\Sigma $ by the definition of induced contact type and by assumption 
$[p_idq_i-\lambda ]=0\in H^1(\Sigma )$, we know that 
\begin{equation}
p_idq_i=\lambda +d\beta (\sigma ) \ on \ \Sigma 
\end{equation}
Then by the proof of Lemma3.3, one computes 
\begin{eqnarray}
G'^*(p_idq_i)&=&G'^*(\lambda )+G'^*d\beta =
\eta (\cdot ,\cdot )^*\lambda +
\eta (\cdot ,\cdot )^*d\beta \nonumber \cr 
&=&\eta ^*_s\lambda +i_X\lambda ds +d\eta ^*\beta 
=ds+d\beta 
\end{eqnarray}
here we also use $\beta $ denote the $\eta ^*\beta$. 

All above construction is contained in \cite{ma}. Now 
we introduce the upshot construction in \cite{mo}:  
\begin{eqnarray}
&&F_0':{\cal {L}}\times R\times R\to (R\times \Sigma )\cr 
&&F_0'(((l,s,a)=(a,\eta _s(l)))
\end{eqnarray}
Now we embed a elliptic curve $E$ long along $s-axis$ and thin along $a-axis$ such that 
$E\subset [-K,K]\times [0,\varepsilon]$. We parametrize the $E$ by $t'\in S^1$.

\begin{Lemma}
If there does not exist any Reeb chord 
in $(\Sigma ,\lambda )$, 
then 
\begin{eqnarray}
&&F_0:{\cal {L}}\times S^1\to (R\times \Sigma )\cr 
&&F_0(l,t')=(a(t'),\eta _{s(t')}(l))
\end{eqnarray}
is a compact Lagrangian submanifold. Moreover 
\begin{eqnarray}
&&l(R\times \Sigma ,F_0({{\cal {L}}}\times S^1),d(e^a\lambda ))\cr 
&&=
\inf \{\int _Df^*de^a\lambda >0|f:(D,\partial D)\to (R\times \Sigma,F_0({\cal {L}}\times S^1)\}\cr 
&&=area(E)
\end{eqnarray}
\end{Lemma}
Proof. We check
that 
\begin{eqnarray}
{F_0}^*(e^a\lambda )&=&e^{a(t')}ds(t')
\end{eqnarray}
So, $F_0$ is a Lagrangian embedding.

If the circle $C$ homotopic to $C_1\subset {\cal {L}}\times s_0$ then  we compute
\begin{eqnarray}
\int _CF_0^*(e^a\lambda )=\int _{C_1}F_0^*(e^a\lambda )=0. 
\end{eqnarray}
since $\lambda |C_1=0$ due to $C_1\subset {\cal {L}}$ and 
$\cal L$ is Legendre submanifold. 

If the circle $C$ homotopic to $C_1\subset l_0\times S^1$ then  we compute
\begin{eqnarray}
\int _CF_0^*(e^a\lambda )=\int _{C_1}F_0^*(e^{a(t')}ds(t'))=n(area(E)). 
\end{eqnarray}
This proves the Lemma.

Now we modify the above construction as follows:  
\begin{eqnarray}
&&F':{\cal {L}}\times R\times R\to 
([0,\varepsilon ]\times \Sigma )\subset T^*M\cr 
&&F'(l,s,a)=(a,\eta _s(l))
\end{eqnarray}
Now we embed a elliptic curve $E$ long along $s-axis$ and thin along $b-axis$ such that 
$E\subset [-s_1,s_2]\times [0,\varepsilon]$. We parametrize the $E$ by $t'$.

\begin{Lemma}
If there does not exist any Reeb chord 
in $(\Sigma ,\lambda )$, 
then 
\begin{eqnarray}
&&F:{\cal {L}}\times S^1\to ([0,\varepsilon ]\times \Sigma )\subset 
T^*M\cr 
&&F(l,t')=(a(t'),\eta _{s(t')}(l))
\end{eqnarray}
is a compact Lagrangian submanifold. Moreover 
\begin{equation}
l(V',F({{\cal {L}}}\times S^1),d(p_idq_i))=area(E)
\end{equation}
\end{Lemma}
Proof. We check
that 
\begin{eqnarray}
F^*(p_idq_i)={F}^*(e^a\lambda +d\beta ).
\end{eqnarray}
This proves the Lemma.

Now we construct an isotopy of Lagrangian embeddings as follows:  
\begin{eqnarray}
&&F':{\cal L}\times S^1\times [0,1]\to V'\cr 
&&F'(l,t',t)=h^T_t(a(t'),\eta _{s(t')}(l))  \cr 
&&F'_t(l,t')=F'(l,t',t).
\end{eqnarray}

\begin{Lemma}
If there does not exist any Reeb chord for $X_\lambda $ 
in $(\Sigma ,\lambda )$ then $F'$ is an weakly exact isotopy of Lagrangian 
embeddings. Moreover 
for the choice of $T=T_\Sigma $ satisfying 
$[0,\varepsilon ]\times \Sigma \cap h^T_1([0,\varepsilon ]\times \Sigma )=\emptyset $, then 
$F'_0({\cal L}\times S^1)\cap F'_1({\cal L}\times S^1)=\emptyset $.
\end{Lemma}
Proof. By Lemma3.1-3.5 and below.

\vskip 5pt

   Let 
$(V',\omega ')=(T^*M, dp_i\wedge dq_i)$, 
$W'=F({\cal L}\times S^1)$, and 
$(V,\omega )=(V'\times C,\omega '\oplus \omega _0)$. 
As in \cite{gro}, we use figure eight trick invented by Gromov to 
construct a Lagrangian submanifold in $V$ through the 
Lagrange isotopy $F'$ in $V'$. 
Fix a positive $\delta <1$ and take a $C^{\infty }$-map $\rho :S^1\to 
[0,1]$, where the circle $S^1$ is parametrized by $\Theta \in [-1,1]$, 
such that the $\delta -$neighborhood $I_0$ of $0\in S^1$ goes to 
$0\in [0,1]$ and $\delta -$neighbourhood $I_1$ of $\pm 1\in S^1$ 
goes $1\in [0,1]$. 
Let $h^T_\rho (t,w')=h^T_{\rho (t)}(w')$ and  
\begin{eqnarray}
\tilde {l}&=&h^{T*}_{\rho }(p_idq_i)
=p_idq_i-\rho (\Theta )Tdg\cr 
&=&e^{a(t')}ds(t')+d\beta )-\rho Tdg=e^{a(t')}ds(t')+d\beta +d\rho Tg)+Tgd\rho \cr 
&=&e^{a(t')}ds(t')+d\beta +d\rho Tg)-Tg\rho '(\Theta )d\Theta \cr
&=&e^{a(t')}ds(t')+d\beta +d\rho Tg)-\Phi d\Theta  
\end{eqnarray}
be the pull-back of the form 
$\tilde {l}'=e^{a(t')}ds(t')+d\beta +d\rho Tg)
-\psi (s,t)dt $ to $W'\times S^1$ under the map 
$(w',\Theta )\to (w',\rho (\Theta ))$ and 
assume without loss of generality $\Phi $
vanishes on $W'\times (I_0\cup I_1)$. 
Since 
$[\tilde {l}']|W'\times \{t\}=[e^{a(t')}ds(t')]$ is independent 
of $t$, so $F'$ is weakly exact. It is crucial here 
$|-\psi (s,t)|\leq M_0$  and $M_0$ is independent of $area(E)$.

  Next, consider a map $\alpha $ of the annulus $S^1\times [\Phi _-,\Phi _+]$ 
into $R^2$, where $\Phi _-$ and $\Phi _+$ are the lower and the upper 
bound of the fuction $\Phi $ correspondingly, such that 
   
   $(i)$ The pull-back under $\alpha $ of the form 
$dx\wedge dy$ on $R^2$ equals $-d\Phi \wedge d\Theta $. 
  
   $(ii)$ The map $\alpha $ is bijective on $I\times [\Phi _-,\Phi _+]$ 
where $I\subset S^1$ is some closed subset, 
such that $I\cup I_0\cup I_1=S^1$; furthermore, the origin 
$0\in R^2$ is a unique double point of the map $\alpha $ on 
$S^1\times 0$, that is 
$$0=\alpha (0,0)=\alpha (\pm 1,0),$$  
and 
$\alpha $ is injective on $S^1=S^1\times 0$ minus $\{ 0,\pm 1\}$. 

   $(iii)$ The curve $S^1_0=\alpha (S^1\times 0)\subset R^2$ ``bounds'' 
zero area in $R^2$, that is $\int _{S^1_0}xdy=0$, for the $1-$form 
$xdy$ on $R^2$. 
\begin{Proposition}
Let $V'$, $W'$ and $F'$ as above. Then there exists  
an exact Lagrangian embedding $F:W'\times S^1\to V'\times R^2$ 
given by $F(w',\Theta )=(F'(w',\rho (\Theta )),\alpha (\Theta ,\Phi ))$.
Denote $F(W'\times S^1)$ by $W$. $W\subset T^*M\times B_{r_0}(0)$ with 
$4\pi r_0^2=8M_0$. 
\end{Proposition}
Proof. Similar to \cite[2.3$B_3'$]{gro}.

\subsection{Formulation of Hilbert manifolds}

Let $(\Sigma ,\lambda )$ be a closed $(2n-1)-$ dimensional manifold
with a contact form $\lambda $ of induced type 
in $T^*M$, it is well-known that $T^*M$ is a Stein manifold, so 
it is exausted by a proper pluri-subharmonic function. In fact 
since $M$ 
is an open manifold one can take a proper Morse 
function $g$ on $M$ and let $f={{|p|^2}\over {2}}+\pi ^*g$. 
Then $f$ is pluri-subharmonic function 
on $T^*M$ for some complex structure $J'$ on $T^*M$ tamed by $dp_i\wedge dq_i$(see\cite{egc}). 
Since $\Sigma $ is compact and $W'=G'({\cal {L}}\times R)$ 
is contanied in $\Sigma $, by our construction 
we have $W'$ is contained in a compact set $f_c$ 
for $c$ large enough. 

\vskip 3pt

Let $V'=T^*M$ and we choose 
an almost complex structure $J'$ on 
$T^*M$ tamed by $\omega' =dp_i\wedge dq_i$ and 
the metric $g'=\omega '(\cdot , J'\cdot )$(see\cite{gro}). 
By above discussion we know that all mechanism such as 
$W'$ or $\Sigma $ contained in $f_c=\{v'\in T^*M|
f(v')\leq c\}$ for $c$ large 
enough, i.e., contained in a compact set $V'_c$ in $T^*M$. 
Then we expanding near $\partial f^{-1}(c)$ to get 
a complete exact symplectic manifold with 
a complete Riemann metric with injective radius 
$r_0>0$(see\cite{ma}).

\vskip 3pt 

In the following we denote by 
$(V,\omega )=(V'\times R^2,d(p_idq_i )
\oplus dx\wedge dy))$ 
with the metric $g=g'\oplus g_0 $ induced by 
$\omega (\cdot ,J\cdot )$($J=J'\oplus i$ and 
$W\subset V$ a Lagrangian submanifold which was constructed in 
section 3.1, moreover we can slightly perturb the $J'\oplus i$ 
near $p$ such that $J\oplus i$ is integrable near $p$.

   Let 
$${\cal D}^k(V,W,p)=\{ u \in H^k(D,V)|
u(x)\in W \ a.e \ for \ x\in \partial D \ and \ u(1)=p\}$$
for $k\geq 100$.
\begin{Lemma}
Let $W$ be a Lagrangian submanifold in 
$V$. Then, 
$${\cal D}^k(V,W,p)=\{ u \in H^k(D,V)|
u(x)\in W \ a.e \ for \ x\in \partial D \ and \ u(1)=p\}$$
is a pseudo-Hilbert manifold with the tangent bundle
\begin{equation}
T{\cal D}^k(V,W,p)=\bigcup _{u\in {\cal {D}}^k(V,W,p)}
\Lambda ^{k-1}(u^*TV,u|_{\partial D}^*TW,p)
\end{equation}
here 
$$\Lambda ^{k-1}(u^*TV,u|_{\partial D}^*TW,p)=$$ 
$$\{ H^{k-1}-sections \ of \ (u^*(TV),(u|_{\partial D})^*TL)\ 
which \ vanishes  \ at \ 1\} $$
\end{Lemma}
Proof: See \cite{al,kl}. 

\vskip 3pt

   Now we consider  a section
from ${\cal D}^k(V,W,p)$ to $T{\cal D}^k(V,W, p)$ follows as in 
\cite{al,gro}, i.e., 
let $\bar \partial :{\cal D}^k(V,W,p)\rightarrow T{\cal D}^k(V,W,p)$
be the Cauchy-Riemmann section 
\begin{equation}
\bar \partial u={{\partial u}\over {\partial s}}
+J{{\partial u}\over {\partial t}}  \label{eq:CR}
\end{equation}
for $u\in {\cal D}^k(V,W,p)$.

\begin{Theorem}
The Cauchy-Riemann section $\bar \partial $ defined in (\ref{eq:CR})
is a Fredholm section of Index zero.
\end{Theorem}
Proof. According to the definition of the Fredholm section, 
we need to prove that
$u\in {\cal D}^k(V,W,p)$, the linearization
$D\bar \partial (u)$ of $\bar \partial $ at $u$ is a linear Fredholm
operator.
Note that
\begin{equation}
D\bar \partial (u)=D{\bar \partial _{[u]}}
\end{equation}
where
\begin{equation}
(D\bar \partial _{[u]})v=\frac{\partial v}{\partial s}
+J\frac{\partial v}{\partial t}+A(u)v
\end{equation}
with 
$$v|_{\partial D}\in (u|_{\partial D})^*TW$$
here $A(u)$ is $2n\times 2n$
matrix induced by the torsion of
almost complex structure, see \cite{al,gro} for the computation.

   Observe that the linearization $D\bar \partial (u)$ of 
$\bar \partial $ at $u$ is equivalent to the following Lagrangian 
boundary value problem
\begin{eqnarray}
&&{{\partial v}\over {\partial s}}+J{{\partial v}\over {\partial t}}
+A(u)v=f, \ v\in \Lambda ^k(u^*TV)\cr 
&&v(t)\in T_{u(t)}W, \ \ t\in {\partial D}  \label{eq:Lin}
\end{eqnarray}
One 
can check that (\ref{eq:Lin}) 
defines a linear Fredholm operator. In fact, 
by Proposition 2.2 and Lemma 2.1, since the operator $A(u)$ is a compact, 
we know that the operator $\bar \partial $ is a nonlinear Fredholm operator 
of the index zero.

\begin{Definition}
Let $X$ be a Banach manifold and $P:Y\to X$ the Banach 
vector bundle.
A Fredholm section $F:X\rightarrow Y$ is
proper if $F^{-1}(0)$ is a compact set and is called 
generic if $F$ intersects the zero section transversally, see \cite{al,fhv,
gro}.
\end{Definition}
\begin{Definition}
$deg(F,y)=\sharp \{ F^{-1}(0)\} mod2$ is called the Fredholm
degree of a Fredholm section (see\cite{al,fhv,gro}).
\end{Definition}
\begin{Theorem}
The Fredholm section
$F=\bar \partial : {\cal D}^k(V,W,p)\rightarrow T({\cal D}^k(V,W,p))$
constructed in (\ref{eq:CR}) is proper near $F^{-1}(0)$ and 
$$deg(F,0)=1$$
\end{Theorem}
Proof: We assume that $u:D\mapsto V$ be a $J-$holomorphic disk
with boundary $u(\partial D)\subset W$ and 
by the assumption that $u$ is homotopic to the 
constant map $u_0(D)=p$. Since almost complex
structure ${J}$ tamed by  the symplectic form $\omega $,
by stokes formula, we conclude $u: D\rightarrow
V$ is a constant map. Because $u(1)=p$, We know that
$F^{-1}(0)={p}$ which implies the properness.
Next we show that the linearizatioon $DF(p)$ of $F$ at $p$ is
an isomorphism from $T_p{\cal D}(V,W,p)$ to $E$.
This is equivalent to solve the equations
\begin{eqnarray}
&&{\frac {\partial v}{\partial s}}+J{\frac {\partial v}{\partial t}}
+Av=f\cr 
&&v|_{\partial D}\subset T_pW
\end{eqnarray}
here $J=J(p)=i$ and $A(=0)$ a constant zero matrix. By Lemma 2.1, we know that $DF(p)$ is an isomorphism.
Therefore $deg(F,0)=1$.

\section{Non-properness of a Fredholm section}

In this section we shall construct a non-proper Fredholm section 
$F_1:{\cal D}\rightarrow E$ by perturbing 
the Cauchy-Riemann section as in 
\cite{al,gro}.

\subsection{Anti-holomorphic section}

  Let $(V',\omega ')=(T^*M,\omega _M)$ and  
$(V,\omega )=(V'\times C, \omega '\oplus \omega _0)$, 
and 
$W$ as in section3 and $J=J'\oplus i$, $g=g'\oplus g_0$, 
$g_0$ the standard metric on $C$. 

   Now let $c\in C$ be a non-zero vector. We consider the 
equations
\begin{eqnarray}
v=(v',f):D\to V'\times C \nonumber \\
\bar \partial _{J'}v'=0,\bar \partial f=c\ \ \ \nonumber \\
v|_{\partial D}:\partial D\to W\ \ \ 
\end{eqnarray}
here $v$ homotopic to constant map 
$\{ p\}$ relative to $W$. 
Note that $W\subset V'\times B_{r_0}(0)$. 
\begin{Lemma}
Let $v$ be the solutions of (4.1), then one has 
the following estimates
\begin{eqnarray}
E({v})=
\int _D(g'({{\partial {v'}}\over {\partial x}},
{J'}{{\partial {v'}}\over {\partial x}})
+g'({{\partial {v'}}\over {\partial y}},
{J'}{{\partial {v'}}\over {\partial y}}) \nonumber \\
+g_0({{\partial {f}}\over {\partial x}},
{i}{{\partial {f}}\over {\partial x}})
+g_0({{\partial {f}}\over {\partial y}},
{i}{{\partial {f}}\over {\partial y}}))d\sigma 
\leq 4\pi r_0^2. 
\end{eqnarray}
\end{Lemma}
Proof: Since $v(z)=(v'(z),f(z))$ satisfy (4.1)
and $v(z)=(v'(z),f(z))\in V'\times C$ 
is homotopic to constant map $v_0:D\to \{ p\}\subset W$ 
in $(V,W)$, by the Stokes formula
\begin{equation}
\int _{D}v^*(\omega '\oplus \omega _0)=0
\end{equation}
Note that the metric $g$ is adapted to the symplectic form 
$\omega $ and $J$, i.e., 
\begin{equation}
g=\omega  (\cdot ,J\cdot )
\end{equation}
By the simple algebraic computation, we have 
\begin{equation}
\int _{D}{v}^*\omega  ={{1}\over {4}}
\int _{D^2}(|\partial v|^2 
-|\bar {\partial }v|^2)=0
\end{equation}
and 
\begin{equation}
|\nabla v|^2={{1}\over {2}}(
|\partial v|^2 +|\bar \partial v|^2)
\end{equation}
Then 
\begin{eqnarray}
E(v)&=&\int _{D} |\nabla v|^2 \nonumber \\ 
      &=&\int _{D}\{ {{1}\over {2}}(
|\partial v|^2+|\bar \partial v|^2)\} d\sigma \nonumber \\ 
&=&\pi |c|_{g_0}^2
\end{eqnarray}
By the equations (4.1), 
one get 
\begin{equation}
\bar \partial f=c \ on \ D
\end{equation}
We have 
\begin{equation}
f(z)={{1}\over {2}}c\bar z+h(z)
\end{equation}
here $h(z)$ is a holomorphic function on $D$. Note that  
$f(z)$ is smooth up to the boundary $\partial D$, then, by 
Cauchy integral formula
\begin{eqnarray}
\int _{\partial D}f(z)dz&=&{{1}\over {2}}c\int _{\partial D}
\bar {z}dz+\int _{\partial D}h(z)dz \cr
&=&\pi ic
\end{eqnarray}
So, we have 
\begin{equation}
|c|={{1}\over {\pi}}|\int _{\partial D^2}f(z)dz|
\end{equation}
Therefore, 
\begin{eqnarray}
E(v)&\leq &\pi |c|^2
\leq {{1}\over {\pi }}|\int _{\partial D}f(z)dz|^2      \cr
&\leq &{{1}\over {\pi }}|\int _{\partial D}|f(z)||dz||^2   \cr
&\leq &4\pi |diam(pr_2(W))^2 \cr
&\leq &4\pi r_0^2.
\end{eqnarray}
This finishes the proof of Lemma.

\begin{Proposition}
For $|c|\geq 3r_0$, then the 
equations (4.1)
has no solutions. 
\end{Proposition}
Proof. By (4.11), we have 
\begin{eqnarray} 
|c|&\leq &{{1}\over {\pi }}\int _{\partial D}|f(z)||dz|\cr
&\leq &{{1}\over {\pi }}\int _{\partial D}
diam(pr_2(W))||dz| \cr
&\leq &2r_0
\end{eqnarray}
It follows that $c=3r_0$ can not be obtained by 
any solutions.

\subsection{Modification of section $c$}

Note that the section $c$ is not a section of the 
Hilbert bundle in section 3 since $c$ is not 
tangent to the Lagrangian submanifold $W$, we must modify it as follows:

\vskip 3pt 

  Let $c$ as in section 4.1, we define 
\begin{eqnarray}
c_{\chi ,\delta }(z,v)=\left\{ \begin{array}{ll}
c \ \ \ &\mbox{if\  $|z|\leq 1-2\delta $,}\cr
0 \ \ \ &\mbox{otherwise}
\end{array}
\right. 
\end{eqnarray}
Then by using the cut off function $\varphi _h(z)$ and 
its convolution with section 
$c_{\chi ,\delta }$, we obtain a smooth section 
$c_\delta$ satisfying

\begin{eqnarray}
c_{\delta }(z,v)=\left\{ \begin{array}{ll}
c \ \ \ &\mbox{if\  $|z|\leq 1-3\delta $,}\cr
0 \ \ \ &\mbox{if\  $|z|\geq 1-\delta $.}
\end{array}
\right. 
\end{eqnarray}
for $h$ small enough, for the convolution theory see \cite{hor}.

   Now let $c\in C$ be a non-zero vector and 
$c_\delta $ the induced anti-holomorphic section. We consider the 
equations
\begin{eqnarray}
v=(v',f):D\to V'\times C \nonumber \\
\bar \partial _{J'}v'=0,\bar \partial f=c_\delta \ \ \ \nonumber \\
v|_{\partial D}:\partial D\to W\ \ \  \label{eq:4.16}
\end{eqnarray}
Note that $W\subset V\times B_{r_0}(0)$. 
Then by repeating the same argument as section 4.1., we obtain 
\begin{Lemma}
Let $v$ be the solutions of (\ref{eq:4.16}) and $\delta $ 
small enough, then one has 
the following estimates
\begin{eqnarray}
E({v})\leq 4\pi r_0^2. 
\end{eqnarray}
\end{Lemma}
and

\begin{Proposition}
For $|c|\geq 3r_0$, then the 
equations (\ref{eq:4.16})
has no solutions. 
\end{Proposition}

\subsection{Modification of $J\oplus i$}

Let $(\Sigma ,\lambda )$ be a closed contact 
manifold
with a contact form $\lambda $ of 
induced type in $T^*M$. 
Let $J_M$ be an almost complex 
structure on $T^*M$ and 
$J_1=J_M\oplus i$ the almost complex structure on 
$T^*M\times R^2$ tamed by 
$\omega '\oplus \omega _0$. 
Let $J_2$ be any almost complex structure on 
$T^*M\times R^2$.

\vskip 3pt 

Now we consider the almost conplex structure 
on the symplectic fibration $D\times V \to D$ which will 
be discussed in detail in section 5.1., see also 
\cite{gro}.

\begin{eqnarray}
J_{\chi ,\delta }(z,v)=\left\{ \begin{array}{ll}
i\oplus J _M\oplus i\ \ \ &\mbox{if\  $|z|\leq 1-2\delta $,}\cr
i\oplus J_2 \ \ \ &\mbox{otherwise}
\end{array}
\right. 
\end{eqnarray}
Then by using the cut off function $\varphi _h(z)$ and 
its convolution with section 
$J_{\chi ,\delta }$, we obtain a smooth section 
$J_\delta$ satisfying

\begin{eqnarray}
J_{\delta }(z,v)=\left\{ \begin{array}{ll}
i\oplus J_M\oplus i \ \ \ &\mbox{if\  $|z|\leq 1-3\delta $,}\cr
i\oplus  {J_2} \ \ \ &\mbox{if\  $|z|\geq 1-\delta $.}
\end{array}
\right. 
\end{eqnarray}
as in section 
4.2.

Then as in section 4.2, one can also reformulation 
of the equations (\ref{eq:4.16}) and get similar 
estimates of Cauchy-Riemann equations, we leave it 
as exercises to reader.

\begin{Theorem}
The Fredholm sections $F_1=\bar \partial +c_\delta 
: {\cal  {D}}^k(V,W,p) \rightarrow T({\cal {D}}^k(V,W,p))$ is not proper 
for $|c|$ large enough.
\end{Theorem}
Proof. See \cite{al,gro}.

\section{$J-$holomorphic section}

Recall that $W\subset T^*M\times B_{r_0}(0)$ as in section 3. 
The Riemann metric $g$ on $M\times R^{2}$ 
induces a metric $g|W$.

   Now let $c\in C$ be a non-zero vector and 
$c_\delta $ the induced anti-holomorphic section. We consider the 
nonlinear inhomogeneous equations (4.16) and 
transform it into $\bar J-$holomorphic map by 
considering its graph as in \cite{al,gro}.

Denote by $Y^{(1)}\to D\times V$ the bundle of homomorphisms $T_s(D)\to
T_v(V)$. If $D$ and $V$ are given the disk and the almost 
K\"ahler manifold, then
we distinguish the subbundle $X^{(1)}\subset Y^{(1)}$ which consists of
complex linear homomorphisms and we denote $\bar X^{(1)}\to D\times V$ the
quotient bundle $Y^{(1)}/X^{(1)}$. Now, we assign to each $C^1$-map $
v:D\to V$ the section $\bar \partial v$ of the bundle $\bar X^{(1)}$ over
the graph $\Gamma _v\subset D\times V$ by composing the differential of $v$
with the quotient homomorphism $Y^{(1)}\to \bar {X}^{(1)}$. If $c_\delta 
:D\times
V\to \bar X$ is a $H^k-$ section we write $\bar \partial v=c_\delta $ 
for the
equation $\bar \partial v=c_\delta |\Gamma _v$.

\begin{Lemma}
(Gromov\cite{gro})There exists a unique almost complex 
structure $J_g$ on $D\times V$(which 
also depends on the given structures in $D$ and in $V$), such that 
the (germs of) $J_\delta-$holomorphic sections $v:D\to D\times V$ are exactly and 
only the solutions 
of the equations $\bar \partial v=c_\delta $. Furthermore, the 
fibres $z\times V\subset D\times V$ are $J_\delta-$holomorphic(
i.e. the subbundles $T(z\times V)\subset T(D\times V)$ are $J_\delta-$complex) 
and the structure 
$J_\delta|z\times V$ equals the original structure on $V=z\times V$.
Moreover $J_\delta $ is tamed by $k\omega _0\oplus \omega $ for 
$k$ large enough which is independent of $\delta $.
\end{Lemma}

\section{Proof of Theorem 1.1}

\begin{Theorem}
There exists a non-constant $J-$holomorphic map $u: (D,\partial D)\to 
(T^*M\times C,W)$ with $E({u})\leq 4\pi r^2_0.$
\end{Theorem}
Proof.  
By Gromov's $C^0-$converngence theorem and the results in 
section 4 shows the solutions of equations (4.16) must 
denegerate to a cusp curves, i.e., we obtain a Sacks-Uhlenbeck's 
bubble, i.e., $J-$holomorphic sphere or disk with boundary 
in $W$, the exactness of the symplectic form on $T^*M\times R^2$ 
rules out 
the possibility of $J-$holomorphic sphere. For the more detail, see the proof   
of Theorem 2.3.B in \cite{gro}.

\vskip 3pt 

{\bf Proof of Theorem 1.1}. By Theorem 6.1, we know that 
\begin{equation}
l(V',F({{\cal {L}}}\times S^1),d(p_idq_i))=area(E)\leq 4\pi r_0^2
\end{equation}
But if $K$ large enough, $area(E)>8\pi r_0^2$. 
This 
implies the assumption that ${\cal {L}}$ has no self-intersection 
point under Reeb flow does not hold.

\end{document}